\documentclass{amsart}
\usepackage[T1]{fontenc}
\usepackage{amsthm,amsfonts,amsmath,mathrsfs,amssymb}
\usepackage{verbatim}

\newcommand{\im}{\mathrm{i}}
\newcommand{\dif}{\mathrm{d}}

\newcommand{\D}{\mathbb{D}}
\newcommand{\abs}[1]{|#1|}
\newcommand{\Abs}[1]{\left|#1\right|}
\newcommand{\norm}[1]{\|#1\|}
\newcommand{\Norm}[1]{\left\|#1\right\|}

\newcommand{\Set}[1]{\left\{#1\right\}}

\newcommand{\C}{\mathbb{C}}

\newcommand{\Hp}{\mathscr{H}}

\newtheorem{theorem}{Theorem}[section]
\newtheorem{lemma}[theorem]{Lemma}

\theoremstyle{definition}

\theoremstyle{remark}

\numberwithin{equation}{section}

\title[Local interpolation for Dirichlet series]{Local interpolation in \\ Hilbert
spaces of
Dirichlet series}

\author{Jan-Fredrik Olsen}
\address{Department of Mathematical Sciences, Norwegian University of
Science and Technology (NTNU), NO-7491 Trondheim, Norway}
\email{janfreol@math.ntnu.no}

\author{Kristian Seip}
\address{Department of Mathematical Sciences, Norwegian University of
Science and Technology (NTNU), NO-7491 Trondheim, Norway}
\email{seip@math.ntnu.no}

\thanks{The authors are
supported by the Research Council of Norway grant 160192/V30.}


\begin{document}

\begin{abstract}
We denote by $\Hp$ the Hilbert space of ordinary Dirichlet series
with square-summable coefficients. The main result is that a bounded
sequence of points in the half-plane $\sigma >1/2$ is an interpolating
sequence for $\Hp$ if and only if it is an interpolating sequence for the
Hardy space $H^2$ of the same half-plane. Similar local results are
obtained for Hilbert spaces of ordinary Dirichlet series that relate to
Bergman and Dirichlet spaces of the half-plane $\sigma >1/2$.
\end{abstract}

\maketitle

\section{Introduction}

The aim of this paper is to shed more light on the function theory
of the Hilbert space $\Hp$ that consists of all Dirichlet series of
the form
\[f(s)=\sum_{n=1}^\infty a_n n^{-s}, \]
($s=\sigma+\im t$ a complex variable) with
\[\norm{f}_{\Hp}^2=\sum_{n=1}^\infty |a_n|^2 < + \infty. \]
We refer to \cite{hls97} and \cite{hls99} for the basics of $\Hp$.
In particular, it should be stressed that $\Hp$ can be viewed as the
Hardy space $H^2$ of the infinite-dimensional polydisk $\D^\infty$.
From this perspective, when studying $\Hp$, we restrict to a
one-dimensional complex variety of $\D^\infty$.
The intricacies of the function theory of $\Hp$ can be seen as
reflecting this situation: We are dealing with a space of functions
that ``really'' live on the infinite-dimensional set $\D^\infty$.

The definition of $\Hp$ shows that it consists of functions analytic
in the half-plane $\sigma>1/2$. The key to understanding the local
boundary behavior in this half-plane is the following embedding
\cite[p. 140]{montgomery94}, \cite[Theorem 4.11]{hls97}:
\begin{equation}
    \int_\theta^{\theta+1} \left|f\hspace{-0.1 cm}\left(\frac{1}{2}+\im t\right)
    \hspace{-0.06 cm}\right|^2\ \dif t \le C
    \norm{f}_\Hp^2\label{emb}
\end{equation}
with $C$ an absolute
constant independent of $\theta$. Of course, without further
precautions, the embedding makes sense only when the Dirichlet
series converges for $\sigma=1/2$. However, once \eqref{emb} has
been established for, say, Dirichlet polynomials, it follows that
$f(s)/s$ is in $H^2$ of the half-plane $\sigma>1/2$ for every $f$ in
$\Hp$. In other words, an implicit consequence of \eqref{emb} is
that $f$ has nontangential limits almost everywhere on the line
$\sigma=1/2$, and we have therefore a well-defined boundary limit
function which we by convention choose to denote by $f(1/2+\im t)$.

Another way of seeing a link between $\Hp$ and $H^2$ of the
half-plane $\sigma>1/2$ is by comparing reproducing kernels. For
$\Hp$, it is immediate that its kernel $k^{\mathscr{H}}_{w}(s)$ at
the point $w$ is a translation of the Riemann zeta-function:
\begin{equation*}
  k^{\mathscr{H}}_{w}(s) = \zeta(s + \bar{w}).
\end{equation*}
Since the zeta-function has a simple pole of residue $1$ at $s=1$,
we have \begin{equation}
k^{\mathscr{H}}_{w}(s)=\frac{1}{s+\bar{w}-1}+h(s+\bar{w})
\label{local}
\end{equation}
with $h$ an entire function. Here the first term is the kernel of
$H^2$ of the half-plane $\sigma>1/2$ (properly normalized), and so
\eqref{local} says that, near the diagonals, the kernels for the two
spaces coincide modulo a bounded term.

In view of these observations, one might suspect that locally
functions in $\Hp$ look like functions in $H^2$. The main result of
this paper may be seen as a way of quantifying this similarity:
Based on \eqref{emb} and \eqref{local}, we will prove that the
interpolating sequences for the two spaces coincide, provided that
we consider only bounded sequences of interpolation points.

The main result and its proof are presented in Section~2. We will
then, in Section~3, briefly indicate that minor modifications of
\eqref{emb} and \eqref{local} yield similar interpolation results
for a scale of Hilbert spaces studied by J. E. McCarthy in
\cite{mccarthy04}. Section~4 contains some concluding remarks. In
particular, we present there some simple observations on the problem
of describing the unbounded interpolating sequences for $\Hp$,
merely to hint at the complexity of the problem.

\section{Local interpolation in $\Hp$}

Let $H$ be a Hilbert space of functions on some set $\Omega$. We
assume point evaluation $f\mapsto f(w)$ is bounded for each $w$ in
$\Omega$ such that $H$ has a reproducing kernel. We denote this
kernel by $k^H_{w}(s)$ and say that a sequence $S=(s_j)_{j=1}^\infty$ of distinct
points $s_j$ from $\Omega$ is an interpolating sequence for $H$ if
$f(s_j) = a_j$ has a solution $f$ in $H$ whenever
$(a_j\norm{k^H_{s_j}}^{-1})_{j=1}^\infty$ is in  $\ell^2$.

We set $\C_{1/2}=\{s=\sigma+\im t:\ \sigma>1/2\}$, and let
$H^2(\C_{1/2})$ denote the classical Hardy space of this half-plane.
This is the space of all functions $f$ analytic in $\C_{1/2}$ with
\[ \norm{f}^2_{H^2}=\sup_{\sigma>1/2} \frac{1}{2\pi}\int_{-\infty}^\infty
\abs{f(\sigma+\im t)}^2\dif t< + \infty. \] We have nontangential
boundary limits on the line $\sigma=1/2$ for almost every $t$ and
may express the square of the norm as
\[ \norm{f}^2_{H^2}=\frac{1}{2\pi}\int_{-\infty}^\infty
\Abs{f\hspace{-0.1 cm}\left(\frac12+\im t\right)\hspace{-0.06 cm}}^2\dif t. \] The reproducing kernel
of $H^2(\C_{1/2})$ at $w$ is
\[ k^{H^2}_w(s)=\frac{1}{s+\bar{w}-1}. \]
Our main result can now be stated as follows.

\begin{theorem}\label{H2}
Suppose $S$ is a bounded sequence of distinct points from
$\C_{1/2}$. Then $S$ is an interpolating sequence for $\Hp$ if and
only if it is an interpolating sequence for $H^2(\C_{1/2})$.
\end{theorem}

Needless to say, now H.~S.~Shapiro and A.~L.~Shields's $H^2$ version
\cite{shapiroshields61} of L. Carleson's classical interpolation
theorem \cite{carleson58} gives a geometric description of the
bounded interpolating sequences for $\Hp$.

One implication is immediate from \eqref{local} and the fact that
$f(s)/s$ is in $H^2(\C_{1/2})$ whenever $f$ is in $\Hp$. Namely,
when we solve the problem $f(s_j)=a_j$ with $f$ in $\Hp$, we
simultaneously solve the problem $F(s)=a_j/s_j$ with $F$ in
$H^2(\C_{1/2})$. Also, since $S$ is bounded,
$(a_j/\norm{k_{s_j}^{H^2}}_{H^2})_{j=1}^\infty$ is in $\ell^2$ if and only if
$(s_j a_j/\norm{k_{s_j}^{\mathscr{H}}}_{\Hp})_{j=1}^\infty$ is in $\ell^2.$

Let us now assume that the bounded sequence $S$ is an interpolating
sequence for $H^2(\C_{1/2})$. We wish to prove that then $S$ is also
an interpolating sequence for $\Hp$. To begin with, we observe that
it suffices to show that the subsequence
  \begin{equation*}
    S_{\epsilon} = \Set{ s_j=\sigma_j+\im t_j\in S: \ \frac{1}{2}<\sigma_j \le
\frac{1}{2}+\epsilon}
  \end{equation*}
  is an interpolating sequence for $\Hp$ for some small $\epsilon$. Indeed, it is
clear that
  $S\setminus S_\epsilon$ is a finite sequence, which we may write as
\[ S\setminus S_\epsilon=(s_j)_{j=1}^N. \]
The finite interpolation problem $f_0(s_j)=a_j$, $j=1,...,N$ can be
solved explicitly as follows. Choose primes $p_1,...,p_N$ (not
necessarily distinct) such that the product
\[ B(s)=\prod_{j=1}^N \left(1-p_j^{s_j-s}\right) \]
has simple zeros at the points $s_1,...,s_N$. If we set
$B_j(s) = B(s)/(1-p_j^{s_j-s})$
then the finite interpolation problem has solution
\[ f_0(s)=\sum_{j=1}^N a_j \frac{B_j(s)}{B_j(s_j)}.\]
 To solve the full interpolation problem $f(s_j)=a_j$, we can now
 solve
  \begin{equation*}
    f_\epsilon(s_j) = \frac{a_j - f_0(s_j)}{B(s_j)}, \ \ s_j\in
    S_\epsilon,
  \end{equation*}
  so that we obtain the final solution $f=Bf_\epsilon+f_0$.
  Clearly,
  \begin{equation*}
   \left(\frac{a_j -
        f_0(s_j)}{B_0(s_j)}\frac{1}{\norm{k^\Hp_{s_j}}_\Hp}\right)_{s_j \in S_{\epsilon}} \in \ell^2
    \iff \left( \frac{s_j}{\norm{k^\Hp_{s_j}}_\Hp}\right)_{s_j \in S_{\epsilon}} \in \ell^2,
  \end{equation*}
so that we have reduced the problem to showing that $S_\epsilon$ is
an interpolating sequence for $\Hp$.

Our reason for making the transition from $S$ to $S_\epsilon$ is
that it will allow us to make use of the fact that
  \begin{equation} \label{limit:carleson}
  \lim_{\epsilon\to 0} \sum_{s_j \in S_{\epsilon}} \left(\sigma_j -
  \frac{1}{2}\right)=0.
  \end{equation}
We note that \eqref{limit:carleson} is just a consequence of the
trivial fact that an interpolating sequence for $H^2(\C_{1/2})$ is a
Blaschke sequence in $\C_{1/2}$. Since $S$ is a bounded sequence,
this means that
\[ \sum_{s_j \in S}\left(\sigma_j-\frac{1}{2}\right)< + \infty. \]

It seems difficult to obtain a direct solution of the interpolation
problem since we do not have a reasonable substitute for Blaschke
products. We will instead argue by duality, using the following
lemma of R.~P.~Boas \cite{young01}, \cite{boas41}.
\begin{lemma} \label{lem:equiv}
Suppose $(f_j)_{j=1}^{\infty}$ is a sequence of unit vectors in a
Hilbert space $H$. Then the moment problem $\langle f, f_j \rangle_H
=c_j$ has a solution $f$ in $H$ for every sequence
$(a_j)_{j=1}^\infty$ in $\ell^2$ if and only if there is a positive
constant $m$ such that
\begin{equation}\label{dual}
  \Norm{\sum_{j} c_j f_j}_H \ge m
  \norm{(c_j)}_{\ell^2} \end{equation}
for every finite sequence of scalars $(c_j)$.
\end{lemma}
Thus we need to prove \eqref{dual} with $H=\Hp$ and
$f_j=k^\Hp_{s_j}/\| k^\Hp_{s_j}\|_\Hp$ for $s_j$ in $S_\epsilon$.

To simplify the writing, we set $k_s=k_s^\Hp$ and suppress the index
in the norm setting $\norm{f}=\norm{f}_\Hp$. Let $T$ be a positive
number such that $|t_j|\le T-1$ for every $s_j=\sigma_j + \im t_j$
in $S$. We start by using \eqref{emb}:
  \begin{equation*}
    \Norm{\sum_{s_j\in S_\epsilon} c_j \frac{k_{s_j}}{\norm{k_{s_j}}}}^2 \geq
    m_T \int_{-T}^T \Abs{\sum_{s_j\in S_\epsilon} c_j \frac{k_{s_j}(\im t +
        \frac{1}{2})}{\norm{k_{s_j}}}}^2 \dif t.
  \end{equation*}
  The trick is now to replace the kernels of $\Hp$ by the kernels of
  $H^2(\C_{1/2})$. We use \eqref{local} and the triangle inequality:
  \begin{equation}
    \begin{split} \label{eq:H2 første}
    &\left( \int_{-T}^T \Abs{\sum_{s_j\in S_\epsilon}c_j \frac{k_{s_j}(\im t +
        \frac{1}{2})}{\norm{k_{s_j}}}}^2 \dif t \right)^{1/2}
  \geq \left( \int_{-T}^T \Abs{\sum_{s_j\in S_\epsilon} \frac{c_j\norm{k_{s_j}}^{-1}
}{\im t + \bar{s}_j -
      \frac{1}{2}}}^2 \dif t\right)^{1/2} \\
  & \qquad \qquad - \left( \int_{-T}^T \Abs{\sum_{s_j\in S_\epsilon} c_j h\left(\im
t+\bar{s}_j+\frac{1}{2}\right)
    \norm{k_{s_j}}^{-1}}^2 \dif t \right)^{1/2} \hspace{-0.4 cm}.
    \end{split}
  \end{equation}
We split the first term on the right into two pieces:
  \begin{equation} \label{eq:H2 andre}
    \int_{-T}^T \Abs{\sum_{s_j\in S_\epsilon} c_j \frac{\norm{k_{s_j}}^{-1}}{\im t +
\bar{s}_j -
      \frac{1}{2}}}^2 \dif t =
  \left(\int_{-\infty}^\infty - \int_{\abs{t}>T}\right)
  \Abs{\sum_{s_j\in S_\epsilon} c_j \frac{\norm{k_{s_j}}^{-1}}{\im t + \bar{s}_j -
      \frac{1}{2}}}^2
  \dif t.
  \end{equation}
The point is now that the first term on the right in \eqref{eq:H2
andre} is just
  \begin{equation*}
2\pi    \Norm{\sum_{s_j\in S_\epsilon} c_j
\frac{k_{s_j}}{\norm{k_{s_j}}}}_{H^2}^2 \hspace{-0.4 cm},
  \end{equation*}
so that by using the hypothesis on $S$ and Lemma~\ref{lem:equiv}, we
arrive at the inequality
  \begin{equation*} 
    \begin{split}
    \frac{1}{m_T}\Norm{\sum_{s_j\in S_\epsilon} c_j
      \frac{k_{s_j}}{\norm{k_{s_j}}}} \geq
    m' \norm{(c_j)}_{\ell^2} &- \left( \int_{-T}^T \Abs{\sum_{s_j\in S_\epsilon} c_j
h\hspace{-0.1 cm}\left(\im t+\bar{s}_j+\frac{1}{2}\right)
    \norm{k_{s_j}}^{-1}}^2 \hspace{-0.11 cm}\dif t \hspace{-0.02 cm}\right)^{1/2}
    \\ & \hspace{-1 cm}- \left(\int_{\abs{t}>T}
  \Abs{\sum_{s_j\in S_\epsilon} c_j \frac{\norm{k_{s_j}}^{-1}}{\im t + \bar{s}_j -
      \frac{1}{2}}}^2
  \hspace{-0.11 cm}\dif t \hspace{-0.02 cm} \right)^{1/2} \hspace{-0.4 cm}.
  \end{split}
  \end{equation*}

  The two terms that are subtracted on the right are easily
  estimated. Indeed, by the Cauchy--Schwarz inequality, the square of the first term
is bounded by
\[ \sum_{s_j \in S_{\epsilon}} \abs{c_j}
    \norm{k_{s_j}}^{-1} \int_{-T}^T
    \Abs{h\left(\im t+\bar{s}_j-\frac12\right)}^2 \dif t \leq
    B \sum_{s_j \in S_{\epsilon}} \norm{k_{s_j}}^{-2} \sum_{s_j \in S_{\epsilon}}
    \abs{c_j}^2, \]
   with $B$ depending only on $h$, $S$, and $T$.
   The second term is treated in a similar way.
   We use again the Cauchy--Schwarz inequality and make a calculation to see that
its square is bounded by
       \[ \sum_{s_j \in S_{\epsilon}} \abs{c_j} \norm{k_{s_j}}^{-1} \int_{\abs{t}>T}
    \Abs{\frac{1}{\bar{s}_j - \frac{1}{2} + \im t}}^2 \dif t \leq
    2 \sum_{s_j \in S_{\epsilon}} \norm{k_{s_j}}^{-2} \sum_{s_j \in S_{\epsilon}}
\abs{c_j}^2. \\
    \]
By \eqref{local}, $\norm{k_{s_j}}^{-2}\le C (\sigma_j-1/2)$ with $C$
depending only on $S$. In view of \eqref{limit:carleson}, we obtain
\eqref{dual} by choosing $\epsilon$ sufficiently small. This
completes the proof of the theorem.

\section{Local interpolation in McCarthy's spaces $\Hp_\alpha$}

In \cite{mccarthy04}, McCarthy studied a wider class of Hilbert
spaces of Dirichlet series. In particular, he found a scale of such
spaces that resembles the familiar scale to which the classical
Dirichlet, Hardy, and Bergman spaces belong. We shall now see that,
as far as local interpolation is concerned and in accordance with
what was done above, these spaces can be linked to their classical
counterparts in the half-plane $\sigma>1/2$.

We declare $\Hp_\alpha$ with $\alpha\le 1$ to be the space of all
Dirichlet series of the form
\begin{equation*}
  f(s) = \sum_{n=1}^\infty a_n n^{-s}
\end{equation*}
such that \[ \norm{f}_{\Hp_\alpha}^2=\sum_{n=1}^\infty
  \abs{a_n}^2 \log^{\alpha}(n+1) < + \infty. \]
We observe that every space $\Hp_\alpha$ consists of functions
analytic in the half-plane $\sigma>1/2$. For each $\alpha$ the basic
estimate \eqref{emb} can now be transformed into a version that
reveals the local boundary behavior of functions in $\Hp_\alpha$. To
this end, let $\dif m$ denote Lebesgue area measure and $Q_\theta$
the half-strip $\sigma>1/2$, $\theta<t<\theta+1$. Then writing
\eqref{emb} as
\[
\int_\theta^{\theta+1} |f(\sigma+\im t)|^2\ \dif t \le C
\sum_{n=1}^\infty |a_n|^2 n^{2\sigma-1}, \] multiplying both sides
by $(\sigma -
  \frac{1}{2})^{-\alpha - 1}$,
  and integrating from $\frac{1}{2}$ to $+\infty$, we get for $\alpha<0$:
  \begin{equation} \label{bergman}
    \int_{Q_{\theta}} \abs{f(s)}^2 \left(
  \sigma-\frac{1}{2}\right)^{-\alpha -1} \dif m(s)
    \leq C \norm{f}_{\Hp_\alpha}^2.
  \end{equation}
 By a similar computation for $0<\alpha\le 1$, we obtain
  \begin{equation} \label{dirichlet}
    \int_{Q_{\theta}} \abs{f'(s)}^2 \left(
  \sigma-\frac{1}{2}\right)^{-\alpha +1} \dif m(s)
    \leq C \norm{f}_{\Hp_\alpha}^2.
  \end{equation}

The natural counterparts to the spaces $\Hp_\alpha$ in the
half-plane $\sigma>1/2$ are therefore the corresponding
Dirichlet-type spaces ($0<\alpha\le 1$) and Bergman spaces
($\alpha<0$). As there is essentially nothing new in this
discussion, we will be brief with the technical details, and some of
them will be left to the reader.

Along with $\Hp_\alpha$, we define the space $D_\alpha(\C_{1/2})$ as
follows. For $\alpha<0$, it is the weighted Bergman space of
functions analytic in $\C_{1/2}$ such  that
\[\norm{f}_{D_\alpha}^2= \frac{1}{\pi}\int_{\C_{1/2}} |f(s)|^2
\left(\sigma-\frac12\right)^{-\alpha-1} \dif m(s)< + \infty.\] For
$0<\alpha<1$ we let $D_\alpha(\C_{1/2})$ be the Dirichlet-type of
space of functions analytic in $\C_{1/2}$ such that $f(\sigma)\to 0$
when $\sigma\to\infty$ and
\[\norm{f}_{D_\alpha}^2= \frac{1}{\pi} \int_{\C_{1/2}} |f'(s)|^2
\left(\sigma-\frac12\right)^{-\alpha+1} \dif m(s)< + \infty.\]
Finally, we declare $D_1(\C_{1/2})$ to be the Dirichlet space in
$\C_{1/2}$, or more precisely, the collection of functions $f$
analytic in $\C_{1/2}$ such that  \[\norm{f}_{D_1}^2= \frac{1}{2\pi}
\int_{-\infty}^\infty \left|f \hspace{-0.1 cm}\left(\frac12+\im
t\right)\hspace{-0.06 cm} \right|^2\frac{\dif t}{1+t^2} + \frac{1}{\pi}\int_{\C_{1/2}}
|f'(s)|^2 \dif m(s) < + \infty.\] The reproducing kernels for
$D_\alpha(\C_{1/2})$ at $w$ is
\[
  k^{D_{\alpha}}_{w} (s)=c_\alpha (\bar{w} + s -
  1)^{\alpha-1}
  \]
when $\alpha<1$, with $c_\alpha=(-\alpha)2^{-\alpha-1}$ for
$\alpha<0$ and $c_\alpha=2^{\alpha-1}(1-\alpha)^{-1}$ for
$0<\alpha<1$. In the limiting case $\alpha=1$, we have
\begin{equation*}
  k^{D_1}_{w}(s) = \frac{3-2\bar{w}}{1 - 2\bar{w}} \frac{3 +2s}{1
    +2s}\left( \log\frac{(1 + 2\bar{w})(1 + 2s)}{2^3}
  + \log \frac{1}{\bar{w} + s - 1}\right).
\end{equation*}
What is essential here is that for $s$ and $w$ in a bounded set we
have
\begin{equation*}
  k^{D_1}_{w} (s) = \log\frac{1}{s+\bar{w}-1}+ \text{a bounded
  term}.
\end{equation*}

Returning to $\mathscr{H}_{\alpha}$, we observe that its reproducing
kernel may be expressed in terms of a weighted zeta-function,
\begin{equation*}
  k^{\mathscr{H}_{\alpha}}_{w}(s) = \sum_{n =1}^\infty
\frac{n^{-s-\bar{w}}}{\log^{\alpha}(n+1)}.
\end{equation*}
As for the classical zeta-function, we have a singularity at $s=1$,
in agreement with the behavior of the kernel of
$D_\alpha(\C_{1/2})$:
\begin{lemma}
  For $\alpha < 1$, we have
  \begin{equation*}
    \sum_{n = 1}^\infty \frac{n^{-s}}{\log^{\alpha}(n+1)}
    = \Gamma(1-\alpha)(s-1)^{\alpha-1}+\mathcal{O}(1)
  \end{equation*}
  as $s \to 1$. We also have
  \begin{equation*}
    \sum_{n =1}^\infty \frac{n^{-s}}{\log(n+1)} = \log \frac{1}{s-1} + \mathcal{O}(1)
  \end{equation*}
as $s \to 1$.
\end{lemma}
\begin{proof}
  The proof is a calculation analogous to the one for the Riemann zeta-function found
  for instance in \cite{ivic03}. To begin with,
  \begin{eqnarray*}
    \sum_{n=1}^\infty \frac{n^{-s}}{\log^{\alpha}(n+1)} &=& \int_1^{\infty}
\frac{x^{-s}}{\log^{\alpha}
    (x+1)}\dif[x] \\
  &=& \int_1^{\infty} \frac{x^{-s-1}[x]}{\log^{\alpha}(x+1)}\left( s +
    \frac{\alpha}{\log (x+1)} \frac{x}{x+1}\right) \dif x.
  \end{eqnarray*}
  The integral
    \begin{equation*}
        \begin{split}
            \int_1^{\infty} \frac{x^{-s-1}[x]}{\log^{\alpha}(x+1)}
            \left( s + \frac{\alpha}{\log (x+1)} \frac{x}{x+1}
            \right) \\
            -
            \frac{x^{-s}}{\log^{\alpha} (x+1)} & \left( s + \frac{\alpha}{\log (x+1)}
            \right) \dif x
        \end{split}
    \end{equation*}
  converges absolutely and defines an analytic function in the right
  half-plane. We may therefore pass from $[x]$ to $x$ in our
  integral, and ignore the factor $x/(x+1)$.
  For a similar reason, we may replace $\log(x+1)$ by
  $\log x$ and if necessary change the lower limit of integration.

  When $\alpha < 1$, we make the following computation:
   \begin{equation*}
    \int_1^{\infty} \frac{x^{-s}}{\log^{\alpha}x}
     \dif x =
     \Gamma(1-\alpha)(s-1)^{\alpha-1}.
  \end{equation*}
 This gives the desired result for $0<\alpha<1$. When $\alpha<0$, we
 find, using the functional equation for the gamma-function, that
  \begin{equation*}
\int_1^{\infty} \frac{x^{-s}}{\log^{\alpha} x}\left( s +
    \frac{\alpha}{\log x}\right) \dif x=\Gamma(1-\alpha)(s-1)^{\alpha-1}
  \end{equation*}
as well.
  In the limiting case $\alpha = 1$,  we find that
  \begin{eqnarray*}
    \int_2^{\infty} \frac{x^{-s}}{\log x} \dif x
        = \log \frac{1}{s-1} + \mathcal{O}(1)
  \end{eqnarray*}
as $s\to 1$.
\end{proof}

We have now found the appropriate analogues of the basic relations
\eqref{emb} and \eqref{local}. The proof of the following extension
of Theorem~\ref{H2} is essentially a plain rewriting of the proof in
Section~2. We trust that the interested reader may check the
details.

\begin{theorem}\label{Dalpha}
Suppose $S$ is a bounded sequence of distinct points from $\C_{1/2}$
and assume $\alpha\le 1$. Then $S$ is an interpolating sequence for
$\Hp_\alpha$ if and only if it is an interpolating sequence for
$D_\alpha(\C_{1/2})$.
\end{theorem}

A small technical remark is in order. One may distinguish between
interpolating sequences as defined above and so-called universal
interpolating sequences, i.e., sequences $(s_j)_{j=1}^\infty$ for
which $f\mapsto (f(s_j)/\norm{k^H_{s_j}}_H)$ maps $H$ both into and
onto $\ell^2$. In the latter case, one then has the Carleson
embedding
\begin{equation} \label{carleson} \sum_{j=1}^\infty |f(s_j)|^2
\norm{k^H_{s_j}}^{-2}_H\le C \norm{f}^2_H
\end{equation} with some positive constant $C$. In the case of
bounded interpolating sequences for $\Hp$, there is no reason to
make a distinction because every bounded interpolating sequence for
$\Hp$ is also a universal interpolating sequence for $\Hp$. The same
holds true for $\Hp_\alpha$ when $\alpha<0$. However, for
$\Hp_\alpha$ with $0<\alpha\le 1$ this is no longer the case
\cite{bishop94}, \cite{marshallsundberg93}, and one should therefore
make a distinction. Still it is plain that Theorem~\ref{Dalpha}
remains valid if one replaces each occurrence of the string ``an
interpolating sequence'' by ``a universal interpolating sequence''.

There exist geometric descriptions of the (universal)
interpolating sequences for all $\alpha\le 1$. For $\alpha<0$,
Beurling-type density theorems were proved in \cite{seip93}.
Descriptions in terms of Carleson measures were found by W. Cohn in
the case $0<\alpha<1$ \cite{cohn93} and independently by C. Bishop
and by D. Marshall and C. Sundberg in the case $\alpha=1$
\cite{marshallsundberg93}. For further information, we refer to the
monograph \cite{seip04}.

\section{Concluding remarks}

It may be noted that Theorem~\ref{H2} gives the first general
sufficient condition for zero sequences of functions in $\Hp$. The
statement is that for each bounded interpolating sequence $S$ for
$H^2(\C_{1/2})$ there is a function $f$ in $\Hp$ vanishing on $S$.
On the other hand, a simple argument\footnote{We thank E. Saksman
for bringing this point to our attention.} shows that this sequence
$S$ cannot be the zero sequence of any function in $\Hp$. Indeed,
the almost periodicity of the function $t\mapsto f(\sigma+\im t)$
along with Rouch\'{e}'s theorem implies that for any zero
$\sigma_0+\im t_0$ of a function $f$ in $\Hp$ and
$0<\epsilon<\sigma$ there is a positive number $T$ such that every
rectangle $|\sigma-\sigma_0|<\epsilon$, $\theta<t<\theta+T$ contains
a zero of $f$. In particular, this means that no function in $\Hp$
has a nonempty finite zero sequence.

We finally make some remarks about unbounded interpolating sequences
for $\Hp$. Take first an arbitrary interpolating sequence located on
the real line. We may assume with no essential loss of generality
that the sequence is $\sigma_j=1/2+2^{-j}$, $j=1,2,3,...$. Then any
sequence $s_j=\sigma_j+\im t_j$ will give us a Carleson embedding of
the form \eqref{carleson}. To see this, we associate with
\[ f(s)= \sum_{n=1}^\infty a_n n^{-s} \]
the Dirichlet series
\[ f^+(s)=\sum_{n=1}^\infty |a_n| n^{-s}. \]
This means that $|f(s_j)|\le f^+(\sigma_j)$, and since
$\norm{f}_\Hp=\norm{f^+}_\Hp$, we obtain the result from the fact
that $(\sigma_j)_{j=1}^\infty$ yields a Carleson embedding.

We may next observe that the sequence $(s_j)_{j=1}^\infty$ can be
split into a finite number of interpolating sequences for $\Hp$.
(The number of sequences depends only on $(\sigma_j)_{j=1}^\infty$.)
This follows from Gerschgorin's circle theorem. Indeed,
Lemma~\ref{dual} says that it is enough to check that the normalized
Grammian
\[ \left(
\frac{k^{\Hp}_{s_j}(s_l)}{\norm{k^{\Hp}_{s_j}}_{\Hp}\norm{k^{\Hp}_{s_l}}_{\Hp}}\right)_{j,l=1}^\infty
\]
is invertible as a map on $\ell^2$. We observe that in our case the
entries of the matrix decay exponentially and monotonically away
from the main diagonal. Thus by splitting $(\sigma_j)_{j=1}^\infty$
into sufficiently sparse subsequences, we obtain the invertibility
from Gerschgorin's criterion.

To illustrate a different point, we construct the following
sequence. For each positive integer $j$ pick points equi-distributed
on the line segment $\sigma=1/2+2^{-j}$, $0\le t \le 1$, i.e.,
choose
\[ s_{j,l}=\frac12+2^{-j}+\im/l, \quad l=1,2,...,j.
\] Then Carleson's theorem along with our Theorem~\ref{H2} shows
that $(s_{j,l})$ is an interpolating sequence for $\Hp$. In
particular, the Carleson embedding \eqref{carleson} holds. Now if we
move the points vertically and far apart, the Carleson embedding may
fail. This is a consequence of the almost periodicity of
$t\mapsto\zeta(\sigma+\im t)$. If we measure the distance between
two points in terms of the angle between the corresponding
reproducing kernels, this almost periodicity implies that points
that are far apart in the hyperbolic sense of the half-plane may be
arbitrarily close in the geometry induced by $\Hp$.

The conclusion is that the nature of the problem changes quite
dramatically when we remove the a priori assumption that the
sequence $S$ be bounded. Any nontrivial results about unbounded
interpolating sequences and Carleson measures for $\Hp$ should yield
interesting information about the global behavior of functions in
$\Hp$.

\bibliographystyle{hep}
\bibliography{bibliotek}

\begin{thebibliography}{McC04}

\bibitem[Bis94]{bishop94}
C.~Bishop,
\newblock Interpolating sequences for the {D}irichlet space and its
  multipliers,
\newblock Preprint, 1994.

\bibitem[BJ41]{boas41}
R.~P. Boas~Jr., \textsl{ A general moment problem},
\newblock Amer. J. Math \textbf{ 63}, 361--370 (1941).

\bibitem[Car58]{carleson58}
L.~Carleson, \textsl{ An interpolation problem for bounded analytic functions},
\newblock Amer. J. Math \textbf{ 80}, 921--930 (1958).

\bibitem[Coh93]{cohn93}
W.~Cohn, \textsl{ Interpolation and multipliers on {B}esov and {S}obolev
  spaces},
\newblock Complex Variables Theory Appl. \textbf{ 22}, 35--45 (1993).

\bibitem[HLS97]{hls97}
H.~Hedenmalm, P.~Lindqvist and K.~Seip, \textsl{ A {H}ilbert space of
  {D}irichlet series and systems of dilated functions in $\text{L}^2(0,1)$},
\newblock Duke Math. J. \textbf{ 86}, 1--37 (1997).

\bibitem[HLS99]{hls99}
H.~Hedenmalm, P.~Lindqvist and K.~Seip, \textsl{ Addendum to "A {H}ilbert space
  of {D}irichlet series and systems of dilated functions in
  $\text{L}^2(0,1)$"},
\newblock Duke Math. J. \textbf{ 99}, 175--178 (1999), {math.FA/9512211}.

\bibitem[Ivi03]{ivic03}
A.~Ivic,
\newblock \textsl{ The {R}iemann {Z}eta-{F}unction. {T}heory and
  {A}pplications},
\newblock Dover Publications Inc., 2003.

\bibitem[McC04]{mccarthy04}
J.~E. McCarthy, \textsl{ Hilbert spaces of {D}irichlet series and their
  multipliers},
\newblock Trans. Amer. Math. Soc. \textbf{ 356}(3), 881--893 (2004).

\bibitem[Mon94]{montgomery94}
H.~L. Montgomery,
\newblock \textsl{ Ten {L}ectures on the {I}nterface {B}etween {A}nalytic
  {N}umber {T}heory and {H}armonic {A}nalysis}, volume~84 of \textsl{ CBMS
  Regional Conference Series in Mathematics},
\newblock AMS, 1994.

\bibitem[MS93]{marshallsundberg93}
D.~E. Marshall and C.~Sundberg,
\newblock Interpolating sequences for the multipliers of the {D}irichlet space,
\newblock Preprint. Availiable at
  http://www.math.washington.edu/\~{}marshall/preprints/preprints.html, 1993.

\bibitem[Sei93]{seip93}
K.~Seip, \textsl{ Beurling type density theorems in the unit disk},
\newblock Invent. Math. \textbf{ 113}, 21--39 (1993).

\bibitem[Sei04]{seip04}
K.~Seip,
\newblock \textsl{ Interpolation and {S}ampling in {S}paces of {A}nalytic
  {F}unctions}, volume~33 of \textsl{ University Lecture Series},
\newblock American Mathematical Society, Providence, R. I., 2004.

\bibitem[SS61]{shapiroshields61}
H.~S. Shapiro and A.~L. Shields, \textsl{ On some interpolation problems for
  analytic functions},
\newblock Amer. J. Math \textbf{ 83}, 513--532 (1961).

\bibitem[You01]{young01}
R.~M. Young,
\newblock \textsl{ An {I}ntroduction to {N}onharmonic {F}ourier {S}eries},
\newblock Academic Press, New York, {R}evised {F}irst {E}dition edition, 2001.

\end{thebibliography}

\end{document}